\documentclass{article}

\usepackage{amssymb}
\usepackage{amsmath}
\usepackage{graphicx}
\usepackage{amsthm}

\newcommand{\vanden}{\boldsymbol{n}}
\newcommand{\hatvap}{\hat{\boldsymbol{p}}}
\DeclareMathOperator {\E}{E}
\newcommand{\nden}{n}
\newcommand{\ndcero}{\nden_0}
\newcommand{\nncero}{\nnum_0}
\newcommand{\nnum}{N}

\newcommand{\anmae}{\alpha}
\newcommand{\funo}{\mu_1}
\newcommand{\fdos}{\mu_2}
\newcommand{\diff}{\mathrm d}
\newcommand{\nmae}{\varepsilon}
\newcommand{\vannum}{\boldsymbol{N}}

\begin{document}

\title{Estimation of a Probability with Guaranteed Normalized Mean Absolute Error}
\author{Luis Mendo\thanks{E.T.S. Ingenieros de Telecomunicaci\'on, Polytechnic University of Madrid,
28040 Madrid, Spain. E-mail: lmendo@grc.ssr.upm.es. Telephone: +34 91 549 5700. Fax: +34 91 336 7350.}}
\date{September 2009}

\maketitle

\begin{abstract}
The estimation of a probability $p$ from repeated Bernoulli trials is considered in this letter. A sequential approach is followed, using a simple stopping rule. A closed-form expression and an upper bound are obtained for the mean absolute error of the unbiased estimator of $p$. The results given permit the estimation of an arbitrary probability with a prescribed level of normalized mean absolute error.
\end{abstract}


\section{Introduction}
\label{parte: intro}

Sequential estimation of a probability from a set of observations is considered in this letter. This problem arises, among other fields, in Monte Carlo simulation of communication systems, in which performance is usually measured by a bit error rate (BER) or block error rate (BLER). Contrary to fixed-size Monte Carlo methods, in sequential estimation the sample size is not fixed in advance, but is (randomly) determined by the outcome of the simulation, using a certain stopping rule.

In this letter, the observations are assumed to be a sequence of independent Bernoulli trials with probability of success $p$, which is to be estimated. A simple stopping rule known as inverse binomial (or negative binomial) sampling is considered. This rule consists in observing the sequence until a given number $\nnum$ of successes is reached. The resulting number of trials is denoted as $\vanden$. (Random variables are displayed in bold type throughout the letter.) The uniformly minimum variance unbiased estimator of $p$, for $\nnum \geq 2$, is \cite{Prasad82}
\begin{equation}
\label{eq: hat p N-1 n-1}
\hatvap = \frac{\nnum-1}{\vanden-1}.
\end{equation}
For $\nnum \geq 3$ the mean square error (MSE) of \eqref{eq: hat p N-1 n-1} is known to satisfy $\E[(\hatvap-p)^2]/p^2 < 1/(\nnum-2)$ irrespective of $p$  \cite{Prasad82}. Recent works \cite{Mendo06} \cite{Mendo08a} have shown that, for the modified estimator $\hatvap = (\nnum-1)/\vanden$, the confidence level associated with a relative interval of the form $[p/\fdos, p\funo]$ also satisfies a lower bound irrespective of $p$, for $\nnum \geq 3$ and $\funo, \fdos$ not smaller than certain values. The same result holds for the estimator \eqref{eq: hat p N-1 n-1}, albeit for a reduced range of $\funo, \fdos$ values \cite{Mendo08b}. The referred bound can be improved by allowing estimators of the form $\hatvap = \Omega/(\vanden+d)$, where $\Omega$ and $d$ are selectable parameters \cite{Mendo09_ac}.

This letter analyzes the mean absolute error (MAE) of the estimator \eqref{eq: hat p N-1 n-1}, for $\nnum \geq 2$. Compared to the MSE, the MAE is a more natural error measure, and has several advantages \cite{Blyth80} \cite{Gorard05}. It is simpler, it has a clearer meaning, and it is less sensitive to outlying values. Apparently, its lack of use is in large part motivated by the analytical difficulty associated with the absolute value \cite{Blyth80} \cite{Gorard05} \cite{Bar-Lev99}.

\section{Result}
\label{parte: result}

Let $\ndcero$ and $\anmae_\nnum$ be defined as
\begin{align}
\label{eq: ndcero}
\ndcero & = \left\lfloor \frac{\nnum-1}{p} \right\rfloor + 1, \\
\anmae_\nnum & = \frac{2 (\nnum-1)^{\nnum-2} e^{-\nnum+1}} {(\nnum-2)!}.
\end{align}
For $\nnum \geq 2$, the MAE of \eqref{eq: hat p N-1 n-1} satisfies the following.
\begin{align}
\label{eq: MAE igual}
\frac{\E(|\hatvap - p|)}{p} & = 2 \binom{\ndcero-1}{\nnum-1} p^{\nnum-1} (1-p)^{\ndcero-\nnum+1}, \\
\label{eq: MAE lim}
\lim_{p \rightarrow 0} \frac{\E(|\hatvap - p|)} p & = \anmae_\nnum.
\intertext{Furthermore, $\E(|\hatvap - p|) / p$ is a monotonically decreasing function of $p$, and}
\label{eq: MAE menor}
\frac{\E(|\hatvap - p|)}{p} & < \anmae_\nnum \quad \text{for all } p \in (0,1).
\end{align}

\begin{proof}
Let $\E(|\hatvap - p|)/p$ be denoted as $\nmae(p)$. Given $\nnum$, the probability (mass) function of $\vanden$, $f_\nnum(\nden) = \Pr[\vanden = \nden]$, is
\begin{equation}
\label{eq: f}
f_\nnum(\nden)
= \binom{\nden-1}{\nnum-1} p^\nnum (1-p)^{\nden-\nnum}, \quad \nden \geq \nnum.
\end{equation}
The corresponding distribution function is denoted as $F_{\nnum}(\nden)$.

Using the identities
\begin{align}
\frac{f_\nnum(\nden)}{\nden-1} & = \frac{p f_{\nnum-1}(\nden-1)}{\nnum-1}, \\
\E\left[ \frac{1}{\vanden-1} \right] & = \frac{p}{\nnum-1},
\end{align}
the MAE is computed as
\begin{equation}
\label{eq: MAE 1}
\begin{split}
\E(|\hatvap-p|)
& = \sum_{\nden=\nnum}^{\ndcero} f_\nnum(\nden) \left( \frac{\nnum-1}{\nden-1}-p \right) \\
& \quad - \sum_{\nden={\ndcero}+1}^\infty f_\nnum(\nden) \left( \frac{\nnum-1}{\nden-1}-p \right) \\
& = (\nnum-1) \left( 2 \sum_{\nden=\nnum}^{\ndcero} \frac{f_\nnum(\nden)}{\nden-1} - \sum_{\nden={\nnum}}^{\infty} \frac{f_\nnum(\nden)}{\nden-1}  \right) \\
& \quad - p \left( 2 \sum_{\nden=\nnum}^{\ndcero} f_\nnum(\nden) - \sum_{\nden={\nnum}}^{\infty} f_\nnum(\nden) \right) \\
& = 2p[F_{\nnum-1}(\ndcero-1)-F_{\nnum}(\ndcero)].
\end{split}
\end{equation}
Let $b_{\nden,p}(i)$ denote the binomial probability function with parameters $\nden$, $p$ evaluated at $i$. Taking into account that
\begin{equation}
F_{\nnum-1}(\ndcero-1) = F_{\nnum}(\ndcero) + (1-p) b_{\ndcero-1,p}(\nnum-1),
\end{equation}
from \eqref{eq: MAE 1} it is seen that
\begin{equation}
\nmae(p) = 2(1-p)b_{\ndcero-1,p}(\nnum-1),
\end{equation}
which establishes \eqref{eq: MAE igual}. The limit result \eqref{eq: MAE lim} follows from Poisson's theorem \cite[p.~113]{Papoulis02}.

Let $S = \{(\nnum-1)/k, k = \nnum,\nnum+1, \ldots\}$ and $T = (0,1) \setminus S$. For $p \in T$, $\ndcero$ does not change if $p$ is altered by a sufficiently small amount, which implies that $\nmae$ is continuous and differentiable, with
\begin{equation}
\label{eq: der norm MAE}
\frac{\diff \nmae(p)}{\diff p} = 2 \binom{\ndcero-1}{\nnum-1} p^{\nnum-2} (1-p)^{\ndcero-\nnum} (\nnum-1-\ndcero p).
\end{equation}
Substituting \eqref{eq: ndcero}, this expression is seen to be negative. Let $p \in S$, i.e.~$p = (\nnum-1)/k$ for some $k = \nnum,\nnum+1,\ldots$. Although $\ndcero$ has a jump discontinuity at every point of $S$, the function $\nmae$ is continuous, because
\begin{multline}
\lim_{h \rightarrow 0-} \nmae\left( \frac{\nnum-1}{k} + h \right) = \lim_{h \rightarrow 0+} \nmae\left( \frac{\nnum-1}{k} + h \right) \\
= 2 \binom{k-1}{\nnum-1} \frac{(\nnum-1)^{\nnum-1} (k-\nnum+1)^{k-\nnum+1}}{k^k}.
\end{multline}
In addition, $\nmae$ has left and right derivatives at $p = (\nnum-1)/k$; these are given by \eqref{eq: der norm MAE} replacing $\ndcero$ by $k+1$ or $k$ respectively, with the result that the left derivative is negative and the right derivative is $0$. The function $\nmae$ is thus continuous with negative derivative on $T$ and nonnegative one-sided derivatives on $S$. This implies that $\nmae$ is monotonically decreasing, and \eqref{eq: MAE menor} follows.

It is interesting to note that the monotonicity of $\nmae(p)$ and the result \eqref{eq: MAE menor} for $p \in S$ can also be established using a similar procedure to that in \cite{Mendo08a}. For these values of $p$, \eqref{eq: ndcero} simplifies to $\ndcero = (\nnum-1)/p + 1$. Defining
\begin{equation}
\label{eq: x 1}
x = \frac 1 p \ln \frac{\anmae_\nnum} {\nmae(p)},
\end{equation}
the inequality \eqref{eq: MAE menor} is equivalent to $x>0$. It follows that
\begin{equation}
\label{eq: x 2}
\begin{split}
x & = -\frac 1 p \sum_{i=1}^{\nnum-2} \ln \left( 1 - \frac{ip}{\nnum-1} \right) \\
& \quad - \frac 1 p \left( \frac{\nnum-1}{p} - \nnum + 2 \right)
\ln(1-p) - \frac{\nnum-1}{p}.
\end{split}
\end{equation}
The variable $x$ can be written \cite{Mendo08a} as $\sum_{j=0}^\infty x_j p^j$ with
\begin{equation}
\label{eq: x j}
x_j = \frac{1}{(j+1)(\nnum-1)^{j+1}} \sum_{i=1}^{\nnum-2}i^{j+1} + \frac{\nnum-1}{j+2} - \frac{\nnum-2}{j+1}.
\end{equation}
For $\nnum=2$, $x_j$ reduces to $1/(j+2)$, and is thus positive. For $\nnum \geq 3$, substituting the inequality
\begin{equation}
\sum_{i=1}^{\nnum-2}i^{j+1} > \frac{(\nnum-2)^{j+2}}{j+2}
\end{equation}
into \eqref{eq: x j} gives
\begin{equation}
\label{eq: x j mayor}
(j+1)(j+2)x_j > (\nnum-2) \left(1-\frac{1}{\nnum-1} \right)^{j+1} + j-\nnum+3.
\end{equation}
Let $y_j$ denote the right-hand side of \eqref{eq: x j mayor}. Computing $\partial y_j / \partial j$ as if $j$ were a continuous variable, it is seen that $\partial y_j / \partial j>0$. Thus $y_j > y_0 = 1/(\nnum-1) > 0$ for any $j \geq 1$, which implies that all the coefficients $x_j$ are positive. Therefore $x>0$ for $\nnum \geq 2$, and $xp = \sum_{j=0}^\infty x_j p^{j+1}$ is increasing on $S$, from which $\nmae(p) = \anmae_\nnum e^{-xp}$ is decreasing on $S$.
\end{proof}

\section{Discussion}

The result above allows the estimation of a probability $p$ with a prescribed value of the normalized MAE, $\E(|\hatvap-p|)/p$. This value is guaranteed irrespective of the unknown $p$. For example, if a normalized MAE not exceeding $10\%$ is desired, $N=65$ suffices, according to \eqref{eq: MAE menor}.

\begin{figure}%
\centering%
\includegraphics[width = \columnwidth]{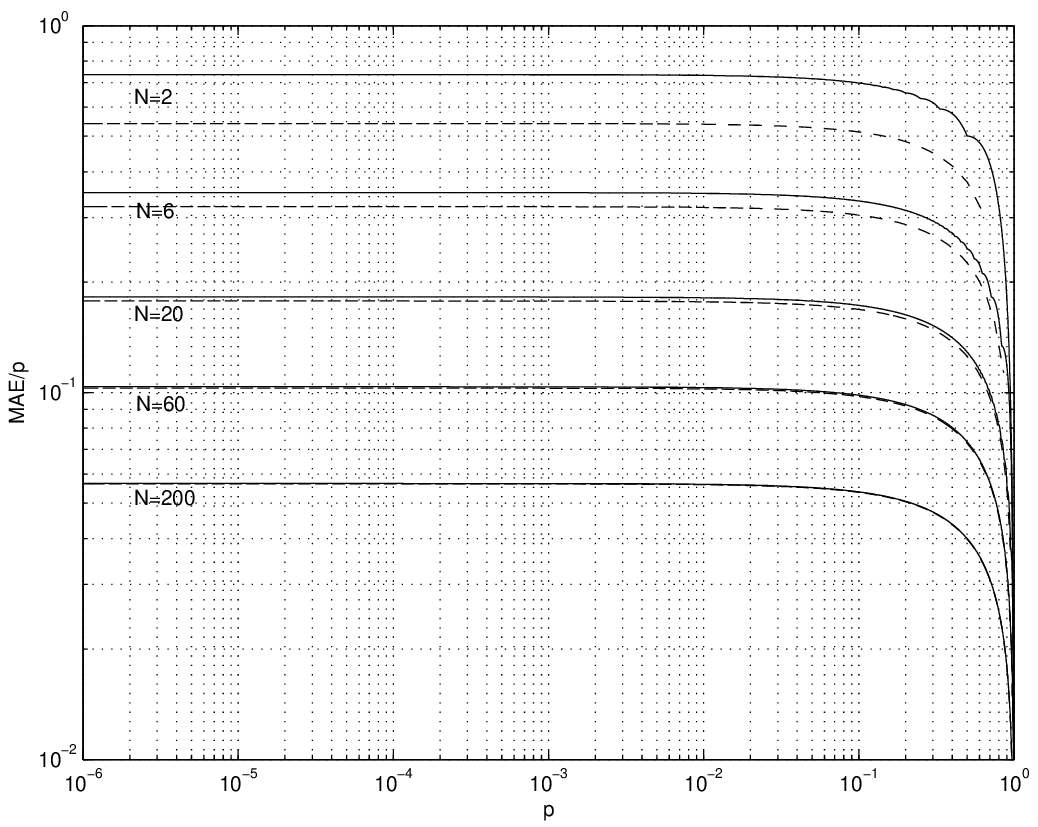}%
\caption{%
\label{fig: MAE_p_fs}%
Normalized MAE as a function of $p$}%
\end{figure}

The behaviour of the normalized MAE as a function of $p$ is depicted in Figure \ref{fig: MAE_p_fs}, with solid lines. The curves show the decreasing character of the normalized MAE. Its nondifferentiability at the points $p=(\nnum-1)/k$, $k \in \mathbb N$ (see proof of the result in Section \ref{parte: result}) can also be clearly observed, specially for low $\nnum$ and large $p$.

\begin{figure}%
\centering%
\includegraphics[width = \columnwidth]{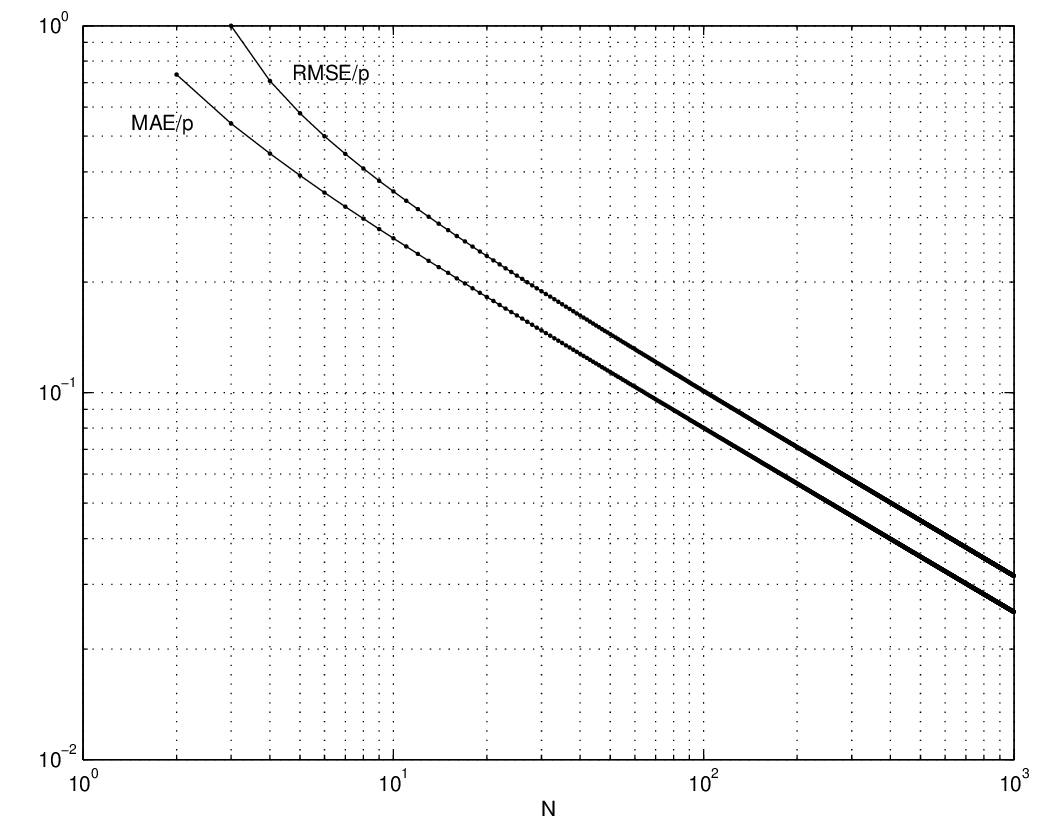}%
\caption{%
\label{fig: bound}%
Bounds on normalized MAE and RMSE as a function of $\nnum$}%
\end{figure}

Figure \ref{fig: bound} shows the bound $\anmae_\nnum$ as a function of $\nnum$. The bound for the root mean square error (RMSE) normalized by $p$, i.e.
\begin{equation}
\frac{\sqrt{\E[(\hatvap-p)^2]}}{p} < \frac{1}{\sqrt{\nnum-2}}
\end{equation}
is also shown for comparison. Both error measures are seen to have the same type of behaviour, with MAE lower than RMSE.

It is interesting to compare \eqref{eq: MAE igual} with the normalized MAE resulting from a fixed sample size $\nden$. In this case, denoting by $\vannum$ the random number of successes, the unbiased estimator $\hatvap = \vannum/\nden$ has a normalized MAE given by \cite[eq.~(1.1)]{Diaconis91}
\begin{equation}
\label{eq: MAE tam fijo}
\frac{\E[|\hatvap-p|]}{p} = 2 \binom{\nden-1}{\nncero-1} p^{\nncero-1} (1-p)^{\nden-\nncero+1}
\end{equation}
with $\nncero = \lfloor \nden p \rfloor + 1$. Since the average sample size in inverse binomial sampling is $\nnum/p$, the comparison is restricted to probabilities $p$ such that $\nnum/p$ is an integer value, and the sample size $\nden$ in the fixed case is taken equal to this value. The resulting fixed-size normalized MAE is shown in Figure \ref{fig: MAE_p_fs} with dashed lines. Dividing \eqref{eq: MAE igual} by \eqref{eq: MAE tam fijo} with $\nden=\nnum/p$, it is easily seen that, for $p \rightarrow 0$, the MAE with inverse binomial sampling is asymptotically $(1+1/(\nnum-1))^{-\nnum+1}e$ times larger than the MAE with fixed sample size. This value is close to $1$ except for very small values of $\nnum$. This is observed in Figure \ref{fig: MAE_p_fs}, which also shows that the MAE ratio is approximately maintained for all values of $p$. It is thus concluded that, in order to guarantee a given normalized MAE, inverse binomial sampling gives an average sample size that is only slightly larger than the sample size that would be necessary in the fixed case (the latter being a function of the unknown $p$).

Possible extensions to this study are: analyzing the effect of imposing a deterministic bound on the number of observations; considering other stopping rules that may be less conservative for $p$ not close to $0$; and replacing the assumption of independent Bernoulli trials by other distributions for the observed variables.

\section{Conclusion}

The MAE in the estimation of a probability $p$ by means of inverse binomial sampling has been characterized. It has been shown that the estimator guarantees a certain value of the normalized MAE irrespective of the unknown $p$. This allows to a priori select a value of the parameter $\nnum$ that meets a prescribed level of normalized error.

The result is quite general, and has many potential applications. In particular, it can be used in simulation and analysis of communication systems, where the performance metric is typically the probability of a certain event.

\section{Acknowledgment}

The author wishes to thank J.~M. Hernando for his valuable help, and the anonymous reviewers for their useful comments.

\end{document}